\documentstyle[12pt, amstex, amssymb, amscd]{amsart}

\newtheorem{theorem}{Theorem}[section]
\newtheorem{vermutung}[theorem]{Vermutung}
\newtheorem{proposition}[theorem]{Proposition}
\newtheorem{korollar}[theorem]{Korollar}
\newtheorem{corollary}[theorem]{Corollary}
\newtheorem{lemma}[theorem]{Lemma}

\theoremstyle{definition}
\newtheorem{definition}[theorem]{Definition}

\theoremstyle{remark}
\newtheorem{remark}[theorem]{Remark}
\newtheorem{remarks}[theorem]{Remarks}
\newtheorem{Bemerkung}[theorem]{Bemerkung}
\newtheorem{Bemerkungen}[theorem]{Bemerkungen}

\def\bt{\begin{theorem}}
\def\et{\end{theorem}}
\def\bl{\begin{lemma}}
\def\el{\end{lemma}}
\def\bp{\begin{proposition}}
\def\ep{\end{proposition}}
\def\bd{\begin{definition}}
\def\ed{\end{definition}}
\def\bv{\begin{vermutung}}
\def\ev{\end{vermutung}}

\newcommand{\BTheorem}{\begin{theorem}}
\newcommand{\ETheorem}{\end{theorem}}
\newcommand{\BVermutung}{\begin{Vermutung}}
\newcommand{\EVermutung}{\end{Vermutung}}
\newcommand{\BFolgerung}{\begin{Folgerung}}
\newcommand{\BBemerkung}{\begin{Bemerkung}}
\newcommand{\EBemerkung}{\end{Bemerkung}}
\newcommand{\BBemerkungen}{\begin{Bemerkungen}}
\newcommand{\EBemerkungen}{\end{Bemerkungen}}

\newcommand{\EFolgerung}{\end{Folgerung}}
\newcommand{\BConjecture}{\begin{conjecture}}
\newcommand{\EConjecture}{\end{conjecture}}
\newcommand{\BProposition}{\begin{proposition}}
\newcommand{\EProposition}{\end{proposition}}
\newcommand{\BCorollary}{\begin{corollary}}
\newcommand{\ECorollary}{\end{corollary}}
\newcommand{\BKorollar}{\begin{korollar}}
\newcommand{\EKorollar}{\end{korollar}}
\newcommand{\BDefinition}{\begin{Definition}}
\newcommand{\EDefinition}{\end{Definition}}
\newcommand{\BLemma}{\begin{lemma}}
\newcommand{\ELemma}{\end{lemma}}
\newcommand{\BRemark}{\begin{remark}}
\newcommand{\BRemarks}{\begin{remarks}}
\newcommand{\ERemark}{\end{remark}}
\newcommand{\ERemarks}{\end{remarks}}

\newcommand{\fg}{{{\frak g}}} 
\newcommand{\fb}{{{\frak b}}}

\newcommand{\fsl}{{{\frak sl}}}

\newcommand{\la}{\lambda}

\newcommand{\al}{\alpha}

\newcommand{\ve}{\varepsilon}

\newcommand{\DZ}{{\Bbb Z}}

\newcommand{\DQ}{{\Bbb Q}}

\newcommand{\Hom}{{\operatorname{Hom}}}

\newcommand{\sur}{\mbox{$\rightarrow\!\!\!\!\!\rightarrow$}}

\numberwithin{equation}{subsection}

\begin{document}

\title{A q-Analogue of Kempf's vanishing theorem}
\author{Steen Ryom-Hansen}
\address{Department of Mathematics\\City University\\
Northampton Square, London, EC1V 0HB, United Kingdom }
\begin{abstract}{We use deep properties of Kashiwara's crystal basis to show that 
the induction functor $ \mbox{H}^0 (-) $ introduced by  Andersen, Polo and Wen satisfies an 
analogon of Kempf's vanishing Theorem.}
\end{abstract}
\email{S.R.H.Ryom-Hansen@@city.ac.uk}
\thanks{Supported
by the TMR-Network Algebraic Lie Theory ERB FMRX-CT97-0100 and 
EPSRC Grant M22536 }
\maketitle 

\section{Introduction}
\subsection{}
Let $ G $ be a reductive connected algebraic group and let $ B $ be a Borel
subgroup. One of central themes of the representation theory of $ G $ is 
the study of the induction functor $ \mbox{H}^0 $ from $ B $ 
representations to 
$ G $ representations. Many of the features of $ \mbox{H}^0 $ in the 
characteristic zero case also hold in the modular case, e.g. the properties 
that $ \mbox{H}^0( \lambda ) \not= 0 $ if and only if $ \lambda \in P^+ $, 
and that the weights of $ \mbox{H}^0( \lambda ) $ are all less than or equal 
to $ \lambda $. On the other hand the Borel-Weil-Bott theorem fails in 
general in the modular case, and hence the simplicity of 
 $ \mbox{H}^0( \lambda ) $ also breaks down in general. Still, we consider 
the $ \mbox{H}^0( \lambda ) $'s to be the fundamental objects of study, the
reason being that their characters, like in the characteristic zero case,
are given by the Weyl character formula. This fact in turn relies on 
the Kempf vanishing theorem, i.e. that
$$  \mbox{H}^i( \lambda ) = 0 \mbox{ for } i > 0 \mbox{ and } \lambda \in
P^+ $$

\subsection{}
In 1979 Andersen and Haboush independently found a short proof of this 
vanishing, see [A] and [H]. Their idea was to show the following isomorphism
$$ \mbox{H}^i(( \, p^r -1)\rho + p^r \, \la ) \cong St_r \otimes 
\mbox{H}^i( \la )^{(r)} $$
where $ St_r $ is a Steinberg module and the superscript denotes the 
$ r $-order Frobenius twist. Because of ampleness properties the left 
hand side is $ O $ for $ r $ sufficiently big; hence 
$ \mbox{H}^i( \la )^{(r)} $ 
must be zero, and thus also $ \mbox{H}^i( \la ) $.

\subsection{}

In [APW 1,2] and [AW] an induction functor $ \mbox{H}^0 $ for quantum 
groups is constructed and studied in great detail. Many of the results 
in these papers rely on specialization to the modular case. In the 
mixed case however, i.e. the case where the ground field is of positive 
characteristic prime to $ l $, these methods fail to give a generalization 
of the Kempf vanishing theorem when $ l < h $, the Coxeter number. And as 
higher ordered Frobenius twists do not exist for quantum groups, also 
the classical method sketched above fails.

\subsection{}

Our approach to the quantum Kempf vanishing theorem is based on some 
properties of the crystal basis proved by Kashiwara in order to obtain 
the refined Demazure character formula in [K2]. In section 2 and 3 we 
discuss these results and in section 4 and 5 we show how they can be applied 
to deduce the Kempf vanishing theorem for the quantum $ \mbox{H}^0 $ in 
all cases.

\subsection{}

It is a pleasure to thank H. H. Andersen for providing crucial input 
to this work. Thanks are also due to G. Lusztig and M. Kashiwara for 
hosting my stays at MIT in the fall 1992 and at RIMS in the spring 1993. 
Finally, it is a great pleasure to thank M. Kaneda for many useful 
discussions.

\subsection{}

Since its appearance and distribution as a preprint in may 1994, there 
has been a couple of developments closely related to and in 
part dependent of our work, we here mention [Ka1], [Ka2] and [W]. As
pointed out by M. Kaneda 
the paper [AP], is incorrect in the generality stated. It is on 
the other hand a key 
reference for certain ampleness properties that are needed in section 4 of 
our work.
This gap was filled in [Ka1, Ka2].
In [Ka2] the universal coefficient theorems are set up that make it 
possible to obtain the vanishing theorem over the ground ring $ {\DZ}
[q,q^{-1}] $. Finally, Woodcock [W] shows
that the idea of using the above mentioned properties of the 
crystal basis to obtain the
quantized Kempf's vanishing theorem also works from a Schur algebra 
point of view.

\section{Notation and some fundamental results}
\subsection{}

Let $ \fg $ be a semisimple complex finite dimensional Lie-algebra. In 
[K1] $ \fg $ is allowed to be a general Kac-Moody algebra, but otherwise 
we shall more or less follow the terminology of that paper. In particular 
$ \{\al_i \}_I $ is the set of simple roots of $ \fg $, $ \{h_i \}_I $ is 
the set of simple coroots, $ P $ is the weight lattice, $ U_q( \fg ) $ is the 
quantized $ \DQ ( q) $-algebra generated by $ e_i, f_i $ where $ i \in I $ 
and $ q^h $ where $ h \in P^* $, $ A $ is the subring of $ \DQ ( q) $
consisting of rational function regular at $ q=0 $, $ V(\la) $ is the 
Weyl module for $ U_q( \fg ) $ with $ v_{\la} $ as a highest weight vector. 

\subsection{}\label{S1}

Assume $ \fg $ has rank one. Then the weight lattice $ P $ is equal to 
$ \DZ $. For $ \la \geq -1 $ the dimension of $ V(\la ) $ is $ \la +1 $:
a basis is $ \{ f_i^{(k)} v_{\la} \, | \, 0 \leq k \leq \la \} $. The action 
of $ U_q( \fg ) $ on this is given by the following formulas:
\begin{equation}
ff^{(k)}v_{\la} = [k+1] f^{(k+1)}v_{\la} 
\end{equation} 
\begin{equation}
e f^{(k)}v_{\la} = [ \la -k + 1 ] f^{(k-1)}v_{\la} 
\end{equation} 
\begin{equation}
q^h f^{(k)}v_{\la} = q^{\la - 2k } f^{(k)}v_{\la}
\end{equation} 
where by convention $ f^{(-1)} v_{\la} = f^{(\la +1 )}  v_{\la} = 0 $. 

\subsection{}

By $ U_q({\fsl}_2 ) $-theory any $ v \in V(\la) $ can be written uniquely 
in the form 
$$ v = \sum_n \, f_i^{(n)} u_n $$
where $ u_n \in V(\la) $ is a weight vector satisfying $ e_i u_n = 0 $. 
Then the operators $ \tilde{e}_i $ and  $ \tilde{f}_i $ on $ V(\la) $ 
are defined in the following way
$$ \tilde{e}_i v := \sum_n f_i^{(n-1)}u_n, \,\,\,\,\,\,\,\,
\,\,\,\,\,\,\,\,
\tilde{f}_i v := \sum_n f_i^{(n+1)}u_n $$
The $ A $-lattice $ L(\la) \subset V(\la) $ is then defined by
$$  L(\la):= A\langle \,  \tilde{f}_{i_1} \tilde{f}_{i_2} \ldots 
\tilde{f}_{i_k} v_{\la} \,|\, i_j \in I, k \geq 0 \,  \rangle.  $$
And $ B(\la) $ is defined as follows 
$$ B(\la) := \pi ( \{\,  \tilde{f}_{i_1} \tilde{f}_{i_2} \ldots 
\tilde{f}_{i_k} v_{\la} \,|\, i_j \in I, k \geq 0 \,  \}) \subset 
 L(\la) / q L(\la) $$
where $ \pi $ is the canonical map $ \pi: L(\la) \rightarrow 
 L(\la) / q L(\la) $. One of the main results of [K1] is then that $ (
L(\la ), B(\la )  ) $ forms a \underline{lower crystal basis}
of $ V(\la) $; this means among other things 
that $ \tilde{e}_i $ and $ \tilde{f}_i $ induce operators on 
$ B(\la) \cup \{ 0 \} $, see theorem 2 of [K1]. 

\subsection{}

The functions $ \ve_i, \varphi_i $ and $ w_i: B(\la) \rightarrow \DZ $ 
are defined in the following way:
$$ \ve_i(b):=\mbox{max}\{n|\tilde{e}_i b \not= 0 \}, \,\,\,\, 
\varphi_i(b):=\mbox{max}\{n|\tilde{f}_i b \not= 0 \}, \,\,\,\,
w_i := \langle \mbox{weight(b)}, h_i \rangle $$
We have the following relations between them
\begin{equation}\label{C1}
w_i(b)=\varphi_i(b) - \ve_i(b) 
\end{equation}
\begin{equation}\label{C2}
\ve_i(\tilde{f}_i b ) = \ve_i(b) +1  \,\,\,\,
\varphi_i(\tilde{f}_i b ) = \varphi_i(b) - 1 
\end{equation}
and $ (L(\la ), B(\la ) ) $ is a \underline{normal} crystal with respect 
to $ \ve_i, \varphi_i $ and $ w_i $.

\subsection{}
The $ \DZ [q,q^{-1}] $-subalgebra of $ U_q( \fg ) $ generated by 
$ f_i^{(n)}, e_i^{(n)} $ and $ q^h, \left\{ \begin{array}{c} q^{h} \\ n 
\end{array}  \right\} $ for $ h \in P^* $ is denoted $ U^{\DZ}_q ( \fg )  $
and the $ U^{\DZ}_q ( \fg )  $-submodule $ V_{\DZ}(\la) $ of $ V(\la) $ 
is by definition $  U^{\DZ}_q ( \fg ) v_{\la} $. Furthermore $ - $ is 
the $ \DQ $-automorphism of $ U^{\DZ}_q ( \fg ) $ given by the formulas 
$$ \overline{e_i} = e_i, \,\,\,\,
\overline{f_i} = f_i, \,\,\,\,
\overline{q^h} = q^{-h}, \,\,\,\, 
\overline{q} = q^{-1}$$
Then Kashiwara has shown, $ (G_l 2 ) $ in section 7.2. of [K1], that 
\begin{equation} \label{G}
V_{\DZ} (\la) \cap L(\la) \cap \overline{L(\la)} \stackrel{\pi}{\cong}
V_{\DZ}(\la) \cap L(\la) / V_{\DZ}(\la) \cap qL(\la) 
\end{equation}
where $ \pi $ is the canonical map. The inverse of $ \pi $ is denoted 
$ G_{\la} $.

\subsection{}

It is known that the crystal $ B(\la) $ is contained in and gives a 
$ \DZ $-basis of $ V_{\DZ}(\la) \cap L(\la) / V_{\DZ}(\la) \cap qL(\la)$. From 
~\ref{G} we now get the following results by applying 
Lemma 7.1.2 of [K1]:
\begin{equation} \label{K1}
V_{\DZ} (\la) \cap L(\la) \cap \overline{L(\la)} \cong \bigoplus_{b \in B(\la)}
\DZ G_{\la} (b) 
\end{equation}
\begin{equation} \label{K2}
V_{\DZ} (\la) \cong \bigoplus_{b \in B(\la)} \DZ[q,q^{-1}] G_{\la} (b)
\end{equation}
\begin{equation} \label{K3}
V (\la) \cong \bigoplus_{b \in B(\la)} \DQ (q) G_{\la} (b)
\end{equation}
\begin{equation} \label{K3}
L(\la) \cong \bigoplus_{b \in B(\la)} A G_{\la} (b) \,\,\,\,\,\,
\overline{L(\la)} \cong \bigoplus_{b \in B(\la)} \overline{A} G_{\la} (b)
\end{equation}

Due to these properties $ \{ G_{\la}(b)|b \in B(\la) \} $ is said to be 
a global basis.

\subsection{}

For $ \Gamma $ a $ \DZ[q,q^{-1}] $-algebra the quantum group $ U_q^{\DZ}(\fg)
\otimes_{\DZ[q,q^{-1}]} \Gamma $ is denoted $ U_{\Gamma}(\fg) $. We use 
the notation $ U_{\Gamma}(\fb) $ and $ U_{\Gamma}({\fb}^-) $ for the 
corresponding Borel subalgebras. 

\subsection{}

Let $ \cal C $ be the [APW] category of integrable 
$ U_q^{\DZ}(\fg) $-modules $ M $. Its objects are the $ U_q^{\DZ}(\fg) $-modules $ M $ having 
a weight space decomposition 
$$ M = \bigoplus_{\la} M_{\la} $$  
where for $\lambda \in P $ the $\lambda$'th weight space $M_{\lambda} $ of $ M $ is 
 
$$ M_{\lambda}=\left\{ \, m \in M \,|\, q^h m = q^{<h, \lambda> } m \,, \, 
 \left\{ \begin{array}{c} q^{h} \\ n 
\end{array}  \right\} m = \left\{ \begin{array}{c} q^{<h, \lambda>} \\ n 
\end{array}  \right\} m \, \mbox{ for all } h \in P^* \right\}$$

\noindent
and such that the following local nilpotency condition holds:
$$ \forall m \in M \mbox{ and } i \in I: \,\,\,\,\,
e_i^{(n)}m = f_i^{(n)}m = 0 \mbox{ for } n >> 0  $$.

Similarly, we consider the category $ {\cal C}^{\geq} $ of 
$ U(\fb) $-modules
that admit a weight space decomposition and satisfy:
$$ \forall m \in M \mbox{ and } i \in I: \,\,\,\,\,
e_i^{(n)}m = 0 \mbox{ for } n >> 0  $$.

The category 
$ {\cal C}^{\leq} $ is defined likewise. 

\medskip

The analogous categories of modules for the specialized quantum groups 
$ U_{\Gamma}(\fg) $,
$ U_{\Gamma}(\fb) $ and $ U_{\Gamma}({\fb}^-) $  
are denoted 
$ {\cal C}_{\Gamma} $,
$ {\cal C}^{\geq}_{\Gamma} $ and $ {\cal C}^{\leq}_{\Gamma} $.

\medskip 
For $ M $ a $ U_q^{\DZ}(\fg) $-module, we define 
$ \operatorname{F}(M) $ 
as the largest submodule of $ M $ which belongs to $ \cal C $. 

\medskip 
If $ M \in \cal C $ we let $ P(M ) $ denote the set of weights of $ M $. 

\subsection{}
Let $ N \in {\cal C}^{\geq} $. Consider the $ \DZ[q,q^{-1}] $-module 
$ \operatorname{Hom}_{ U_q^{\DZ}(\fb)}(\, U_q^{\DZ}(\fg), N \, ):= 
\{ f \in \operatorname{Hom}_{\DZ[q,q^{-1}]}( \, U_q^{\DZ}(\fg), N \,) \, | \,
f(ub) = S^{-1}(b) f(u) \,\, \, \, \, \forall u \in U_q^{\DZ}(\fg), \, \, 
\forall b \in U_q^{\DZ}(\fb) \, \} $ where $ S $ is the antipode map of 
$  U_q^{\DZ}(\fb) $. It has the structure of a $  U_q^{\DZ}(\fg) $-module 
through $ (uf)(m):=f(S(u)m) $. Then the [APW] induction $ \mbox{H}^0(N) $ is 
defined as 
$$ \mbox{H}^0(N):= \mbox{F}(\operatorname{Hom}_{ U_q^{\DZ}(\fb)}
(\, U_q^{\DZ}(\fg), N \, )) $$
(Contrary to what APW do ( and what is the tradition for algebraic groups ) 
we are here inducing from positive Borel groups, this is the reason for 
the difference between our definition and the one in [APW] ). The map 
$ Ev:\mbox{H}^0(N) \rightarrow N; f \mapsto f(1) $ is a $  U_q^{\DZ}(\fb) $-linear map; it induces the Frobenius reciprocity isomorphism:
$$ \operatorname{Hom}_{ U_q^{\DZ}(\fb)}(\, E,N \, ) = 
\operatorname{Hom}_{ U_q^{\DZ}(\fg)}(\, E, \mbox{H}^0(N) \, ) \,\,\,\,\,\,
\forall N \in {\cal C}^{\geq}, E \in \cal C $$
This is the universal property of $ \mbox{H}^0 $. There is also a 
tensor product theorem for $ \mbox{H}^0 $. 

\medskip 
For any $ \DZ[q,q^{-1}]$-algebra $ \Gamma $ we have likewise an induction 
functor $ \mbox{H}^0_{\Gamma}: {\cal C}^{\geq}_{\Gamma} \rightarrow 
{\cal C}_{\Gamma}. $

\section{$W$-filtrations and crystal bases.}
\subsection{}

In this section we shall improve on some of the results of [K2]. In that 
paper all theorems deal with the rings $ \DQ (q) $ and $ \DQ[q,q^{-1}] $;
however we need the results to hold also for $ \DZ[q,q^{-1}] $. 

\subsection{}

Throughout the rest of this section we fix an $ i \in I $ and consider the 
corresponding $ {\fsl}_2 $-component $ {\fg}_i $ of $ \fg $. Let $ W^l ( \la ) $ 
be the sum of all $  U_q({\fg}_i) $-submodules of $ V(\la) $ of dimension greater 
than or equal to $ l $. Furthermore $ W^l( B(\la) ) $ and $ I^l ( B(\la) ) $ are defined 
as the following subsets of $ B(\la) $:
$$  W^l( B(\la) ):= \{ b \in B(\la) \, | \, \ve_i(b) + \varphi_i(b) \geq l \}$$
$$  I^l( B(\la) ):= \{ b \in B(\la) \, | \, \ve_i(b) + \varphi_i(b) =l \}$$
where $ \ve_i $ and $ \varphi_i $ are the functions mentioned in 2.3. Then (3.3.1) of 
[K2] says that 
\begin{equation}
\label{W}
W^l(\la) = \bigoplus_{b \in W^l( B(\la) ) } \DQ (q) G_{\la} (b) 
\end{equation}
We can improve this to the following Lemma:
\begin{lemma}\label{W1} $ W^l(\la) \cap V_{\DZ}(\la) = 
\bigoplus_{b \in W^l( B(\la) ) } 
\DZ [q,q^{-1}] G_{\la} (b) $ 
\end{lemma}
\begin{pf*}{Proof}
The inclusion $ \supset $ follows from ~\ref{W} together with ~\ref{K2}. For 
the other inclusion assume $ w \in  W^l(\la) \cap V_{\DZ}(\la) $. Then 
using ~\ref{W} and ~\ref{K2} once more $ w $ can be written as 
$$ w = \sum_{b \in W^l( B(\la) )} f_b \, G_{\la} (b) 
 = \sum_{b \in  B(\la) )}  g_b \,  G_{\la} (b),  \,\,\,\,\,\,\,\,\,\,
f_b \in \DQ (q),\, g_b \in \DZ [q,q^{-1}] $$
But the $  G_{\la} (b) $'s are independent so we can conclude that 
$ f_b = g_b $. The Lemma is proved.
\end{pf*}

\subsection{}

For $ b \in I^l( B(\la) ) $ we have the following formulas, 3.1.2 of [K2]:
\begin{equation} \label{V1}
f_i^{(k)} G_{\la} (b) \equiv \left[ \begin{array}{c} 
\ve_i(b) + k \\ k  \end{array} \right] \, G_{\la} (\tilde{f}_i^k b)
\,\,\,\,\,\,\,\,\, \mbox{ mod } W^{l+1}(\la) 
\end{equation}
\begin{equation} \label{V2}
e_i^{(k)} G_{\la} (b) \equiv \left[ \begin{array}{c} 
\ve_i(b) + k \\ k  \end{array} \right] \, G_{\la} (\tilde{e}_i^k b)
\,\,\,\,\,\,\,\,\, \mbox{ mod } W^{l+1}(\la) 
\end{equation}
\begin{lemma} \label{A} The above formulas also hold mod $ W^{l+1}(\la) \cap 
V_{\DZ} (\la) $. 
\end{lemma}
\begin{pf*}{Proof}
We have that $ G_{\la}(b) \in V_{\DZ} (\la) \, \, \, 
 \forall b \in B(\la) $ and hence
$ f_i^{(k)} G_{\la} (b) -\left[ \begin{array}{c} 
\ve_i(b) + k \\ k  \end{array} \right] \, G_{\la} (\tilde{f}_i^k b) \in 
V_{\DZ} (\la) $. This proves the Lemma. 
\end{pf*}

\subsection{}
Let us now consider a $ U_q^{\DZ} ( {\fb}_i ) $-module $ N \subset V(\la ) $ 
and a $ U_q^{\DZ} ( {\fb}_i ) $-submodule $ N_{\DZ} $. Assume furthermore 
the existence of a $ B_N \subset B(\la) $ such that 
\begin{equation}
N_{\DZ} \cong \bigoplus_{b \in B_N  } \, \DZ[q,q^{-1}]  G_{\la} (b)
\end{equation}
\begin{equation}
N \cong \bigoplus_{b \in B_N  } \, \DQ (q)  G_{\la} (b)
\end{equation}
According to Lemma 3.1.2 of [K2] $ B_N $ must then satisfy that
\begin{equation} \label{BN}
\tilde{e}_i B_N \subset B_N \sqcup \{0\} 
\end{equation}
We now make the following definitions:
$$ \widetilde{N}_{\DZ} = \sum_{n} f_i^{(n)} N_{\DZ}, \,\,\,\,\,\,\, 
\widetilde{N} = \sum_{n} f_i^{(n)} N, \,\,\,\,\,\,\,
\widetilde{B}_N = \bigcup_{n} \tilde{f}_i^{n} B_N $$
Then Kashiwara has shown, Theorem 3.1.1 of [K2], that 
\begin{equation} \label{N}
\widetilde{N} = \bigoplus_{b \in \widetilde{B}_N  } \, \DQ (q) 
G_{\la} (b)
\end{equation}
We wish to improve this to a statement about $ \widetilde{N}_{\DZ} $. 

\begin{lemma} $  \widetilde{N}_{\DZ} = \bigoplus_{b \in 
\widetilde{B}_N  } \, \DZ[q,q^{-1}]  G_{\la} (b) $ \label{L1}
\end{lemma}
\begin{pf*}{Proof} The intersection of the right hand side of ~\ref{N} 
with $ V(\la)_{\DZ} = \bigoplus_{b \in 
B(\la) } \, \DZ[q,q^{-1}]  G_{\la} (b) $ equals the right hand side of 
the Lemma. Hence we must prove that 
\begin{equation}
\widetilde{N} \cap V_{\DZ}(\la ) = \widetilde{N}_{\DZ}
\end{equation}
The inclusion $ \supset $ is clear. For the other inclusion choose an 
$ n \in \widetilde{N} \cap V_{\DZ}(\la) $. Then we can write $ n $ in the 
following two ways
$$ (*)  \,\,\,\,\,\,\,\,\,\,\,
n = \sum_{k \geq 0,\, b \in B_N} c_{k,b}\, f_i^{(k)}  G_{\la} (b) 
= \sum_{ \beta \in B(\la)} d_{\beta}  \, G_{\la} (\beta) $$
where $ c_{k,b} \in \DQ (q) $ and $ d_{\beta} \in  \DZ[q,q^{-1}] $. 

\medskip
We wish to modify the first sum so that the occuring $ b $'s all satisfy 
$ \tilde{e}_i b = 0 $. Assume that $ b $ occurs in the sum and that 
$ \tilde{e}_i b \not= 0 $. Choose $ l $ minimal such that 
$ b \in W^l( B(\la) ) $. By ~\ref{W} we then have $  G_{\la} (b) \in 
W^l ( \la ) $ and thus $  G_{\la} (b) \in W^l ( \la ) \cap \widetilde{N} = 
W^l (\widetilde{N}) $; so it follows from $ {\fsl}_2 $-theory that 
\begin{equation} \label{last1}
G_{\la} (b) \equiv c f_i e_i G_{\la} (b) \,\,\,\,\,\,\, \mbox{ mod }
W^{l+1} ( \widetilde{N} ) 
\end{equation}
where $ c \in \DQ (q) $. On the other hand ~\ref{W} and ~\ref{V2} give that 
\begin{equation} \label{last}
e_i \,  G_{\la} (b)= c_1 \, G_{\la} (\tilde{e}_i b) + 
\sum_{ b \in W^{l+1}(B (\la))} c_{b} \,  G_{\la} (b) 
\end{equation}
with $ c_1, c_b \in \DQ ( q ) $. As $ e_i \, G_{\la} (b) \in N $, it can also 
be written as a $ \DQ (q ) $-combination of the $  G_{\la} (b) $'s with
$ b \in B_N $; hence we get from the independence of the $ G_{\la} (b)$'s 
that $ \tilde{e}_i \, b $ along with the $ b $'s in the sum belong to $ B_N $.
( This is actually the proof of ~\ref{BN} ). Thus, writing ~\ref{last}
in the form $$ e_i G_{\la} (b) \equiv c_1 G_{\la} (\tilde{e}_i b) 
\,\,\,\,\,\,\,\,\, \mbox{ mod } W^{l+1}( \widetilde{N} ) $$
we deduce that 
$$ c\, f_i e_i G_{\la} (b) \equiv  c\, c_1 \, f_i  G_{\la} (\tilde{e}_i b)
\,\,\,\,\,\,\,\,\, \mbox{ mod } 
W^{l+1}( \widetilde{N} ) $$
Combining this with ~\ref{last1} we obtain 
\begin{equation}
G_{\la} (b) \equiv c c_1 f_i G_{\la} (\tilde{e}_i b ) \,\,\,\,\,\,\,\,\, 
\mbox{ mod } 
W^{l+1}( \widetilde{N} ) 
\end{equation}
Using this and descending induction on $ l$ we can write $ n $ as 
promised $$ n = 
\sum_{ b \in B_N \tilde{e}_i b = 0 } c_{k,b} f_i^{(k)}  G_{\la} (b) $$
Using $ {\fsl}_2$-theory once more and an induction like the previous one 
we can furthermore ensure that the occuring $ f_i^{(k)} G_{\la} (b) $ 
satisfy $ k \leq w_i (b) $. Now from ~\ref{V1} and ~\ref{N} we have that 
$$ f_i^{(k)} G_{\la} (b) \equiv  G_{\la} (\tilde{f}_i^k b)
\,\,\,\,\,\,\,\,\, \mbox{ mod } W^{l+1}(\widetilde{N}) $$
Thus the set of vectors in $ \widetilde{N} $ 
$$ \{ \,  f_i^{(k)} G_{\la} (b) \, | \, k \leq w_i(b), b \in B_N, 
\tilde{e}_i b = 0 \,  \} $$
is a $ \DQ ( q ) $-basis of $ \widetilde{N} $ and the base change matrix
from $ \{ \, G_{\la} (b)\, |\, b \in \widetilde{B}_N \} $ is triangular with ones
on the diagonal with respect to a proper indexing of the basis. But then 
in (*) we must have that $ c_{k,b} \in \DZ [q,q^{-1} ] $ and thus 
$ n \in \widetilde{N}_{\DZ} $. This proves the Lemma.
\end{pf*}

\subsection{} \label{VW}

For $ w \in W $ the $ U_q^{\DZ}(\fb) $-submodule $ V_w^{\DZ} (\la) $ of 
$ V^{\DZ} (\la) $ is defined in the following recursive way
$$ V_1^{\DZ} (\la):= v_{\la}, \,\,\,\,\,\,\,\,\,
V_w^{\DZ} (\la):= U_q^{\DZ} ({\fg}_s )  V_{sw}^{\DZ} (\la)  \,\,\,\,\,\,\,\,
sw <w $$
It is shown in Lemma 3.3.1 of [K2] that $ V_w^{\DZ} (\la) $ also has the 
following description
$$ V_w^{\DZ} (\la) =  U_q^{\DZ} ({\fb} ) \, v_{w \la} $$
where $  v_{w \la} \in V_{\DZ} (\la) $ is defined in the following 
recursive way
$$ v_{1 \la} := v_{\la}, \,\,\,\,\,\,\,\,\,\,
v_{w \la} := f^{(m)}_s v_{sw \la} \,\,\,\,\, sw < w $$
where $ m := \langle {\al}_i, sw \la \rangle $. By the quantum Verma 
relations this is indenpendent of the choice of reduced expression of $ w $ 
and hence also $  V_w^{\DZ} (\la) $ is independent of the reduced expression 
of $ w $. 

\medskip

Applying now the above Lemma to $ N_{\DZ} =  V_{sw}^{\DZ} (\la) $ we 
obtain the existence of a $ B_w ( \la ) \subset B(\la ) $ such that 
\begin{equation} \label{VB}
V_w^{ \DZ } ( \la ) = \bigoplus_{ b \in B_w ( \la ) } \DZ [q,q^{-1}] \, 
G_{\la} (b) 
\end{equation}
Then $ B_w ( \la ) $ has the following properties 
\begin{equation}
\tilde{e}_i  B_w ( \la ) \subset  B_w ( \la ) \sqcup \{0 \} 
\end{equation}
\begin{equation}
 B_w ( \la ) = \bigcup_{k} \tilde{f}_s^k  B_{sw} ( \la )
\end{equation}
The first property is a consequence of ~\ref{BN} and the second one follows 
from Lemma~\ref{L1}. Let $ S $ be an $ i $-string, i.e. a subset of $ B(\la)$ 
of the form 
$$ S = \{\, \tilde{f}_i^k b \,|\, k \geq 0, b \in B ( \la ), \tilde{e}_i b = 0 
\, \} $$
where $ b $ is called the highest weight vector. Then $  B_w ( \la ) $ has 
the following property
\begin{equation} \label{St}
 B_w ( \la ) \cap S \mbox{ is either } S \mbox{ or } b \mbox{ or the empty 
set } 
\end{equation} 
This is rather deep; it is the content of Theorem 3.3.3 of [K2].

\subsection{}

We shall investigate the consequences of these properties for the 
$ V_w^{ \DZ } ( \la ) $:
\begin{lemma} \label{Fi} There is a $ U_q^{ \DZ} ( {\fb}_i ) $-filtration of 
$ V_w^{ \DZ } ( \la ) $
$$ 0 = W^l(\la) \, \cap \, V_w^{ \DZ } ( \la ) \subset 
W^{l-1}(\la) \, \cap \, V_w^{ \DZ } ( \la ) \subset \ldots \subset 
W^0(\la) \, \cap \, V_w^{ \DZ } ( \la ) = V_w^{ \DZ } ( \la ) $$
such that the quotients are direct sums of 
Weyl $ U_q^{ \DZ} ( {\fg}_i ) $-modules restricted to 
$ U_q^{ \DZ} ( {\fb}_i ) $ and of rank one $ U_q^{ \DZ} ( {\fb}_i ) $-modules 
having dominant weights.
\end{lemma}
\begin{pf*}{Proof} We can construct a $ U_q^{ \DZ} ( {\fb}_i ) $-filtration 
of $  V_w^{ \DZ } ( \la ) $ in the following way
$$ 0 = W^l(\la) \, \cap \, V_w^{ \DZ } ( \la ) \subset 
W^{l-1}(\la) \, \cap \, V_w^{ \DZ } ( \la ) \subset \ldots \subset 
W^0(\la) \, \cap \, V_w^{ \DZ } ( \la ) = V_w^{ \DZ } ( \la ) $$
where $ l $ is chosen big enough for the first equality to hold. By 
Lemma ~\ref{W1}, ~\ref{VB} and the definition of $ I^l $ the quotients 
are 
$$  W^k(\la) \, \cap \, V_w^{ \DZ } ( \la ) / 
W^{k+1}(\la) \, \cap \, V_w^{ \DZ } ( \la ) \cong 
\bigoplus_{ b \in I^l(B_w ( \la )) } \DZ [q,q^{-1}] \, 
G_{\la} (b) $$
If $ S $ is an $ i $-string then $ \ve_i ( b ) + \varphi_i ( b ) $ is 
constant on $ S $; this follows from ~\ref{C2}. But then $ I^k(S) $ is
either $ S $ or the empty set and we conclude that $ I^k( B_w (\la) ) $ 
inherits the string property ~\ref{St}. If the intersection is $ S $, the 
formula ~\ref{V2} together with Lemma~\ref{A} and the description of the 
Weyl modules for $ U_q^{ \DZ} ( {\fg}_i ) $ in Section~\ref{S1} show that 
$ \{\, G_{\la} (b) \, , \, b \in S \, \} $ gives rise to an 
$ U_q^{ \DZ} ( {\fg}_i ) $ Weyl module restricted to 
$ U_q^{ \DZ} ( {\fb}_i ) $. If the intersection is a highest weight vector 
$ \{b\} $ then it has the weight $ l > 0 $: $ \ve_i (b ) = 0 $ whence 
$ \varphi_i ( b ) = \varphi_i ( b ) + \ve_i (b ) = l $ and then 
~\ref{C1} gives $ w_i(b) = \varphi_i ( b ) - \ve_i (b ) = l $. The Lemma
is proved.
\end{pf*}

\begin{remark} For $ N $ a $ U_q^{ \DZ} ( {\fb}_i ) $-module the filtration
of it by the $ N \, \cap \, W^l ( \la ) $ is denoted the $ W$-filtration of 
$ N $.
\end{remark}

\subsection{} \label{field}

Let $ k $ be a field of characteristic $ p > 0 $ which is made into a 
$ \DZ [q,q^{-1}] $-algebra by sending $ q $ to an $ l $'th root of unity.
( Later on we shall appeal to results of [AW], hence we should really 
impose the restrictions on $ l $ that occur in that paper. However, in [AP]
it is shown that the Frobenius map of [L] can be imployed to get rid of 
these restrictions ). Then, as the filtration quotients in Lemma~\ref{Fi}
are free $  \DZ [q,q^{-1}] $-modules we get by tensoring a filtration 
of $ V_w^k (\la ) :=  V_w^{\DZ} (\la ) \otimes_{ \DZ [q,q^{-1}]} k $ 
having the same properties as the one of $ V_w^{\DZ} (\la )$. 

\section{Joseph's induction functor}

\subsection{}

In this section we shall compare the induction functor $ \cal D $ of 
Joseph, see [J], with the [APW] functor $ \mbox{H}^0 $. Let $ k $ be as 
in Section~\ref{field}. Then $ \cal D $ is defined in the following 
way
\begin{definition} Let $ N $ be a finite dimensional $ U_k( \fb ) $-module
and let $ U_k \supset U_k( \fb ) $ be a parabolic ( in the sense of [APW] ) 
quantum group. Then $ {\cal D} N $ is 
$$ { \cal D } N := \mbox{D}( U_k \otimes_{ U_k( \fb )} N ) $$ 
where the tensor product has structure as a $  U_k $-module through left 
multiplication and $ \mbox{D} $ is the functor from $ U_k $-modules to 
finite dimensional $ U_k $-modules that takes an $ M $ to the largest 
finite dimensional quotient of $ M $ by a $ U_k $-submodule. 
\end{definition}

\begin{remark} Recall that if $ M $ is an integrable $ U_k $-module then 
( the relevant ) Weyl group acts on the weights of $ M $, see [AW] 
Proposition 1.7. Hence, arguing as in [APW], 1.14, we find that there is 
a unique submodule of $  U_k \otimes_{ U_k( \fb )} N $ such that the 
quotient is of maximal dimension, i.e. $ \mbox{D} $ is well defined. 
\end{remark}

The universal property of $ \cal D $ is given by the following 
Frobenius property 
$$ \mbox{Hom}_{U_k( \fb )} ( \, N, E \, ) = 
 \mbox{Hom}_{U_k} ( \,{\cal D} N, E \, ) $$
where $ E, N $ are finite dimensional $ U_k $, $  U_k( \fb )$-modules. 
The isomorphism is induced by the natural $U_k( \fb )$-map $ \sigma: N 
\mapsto {\cal D} N $. Furthermore, $ \cal D $ satisfies a tensor product 
theorem. 

\subsection{}

For $ M $ a $ U_k $-module we set $ M^{dual}:= \mbox{Hom}_k ( \, M, k \, ) $.
Then $ M^{dual} $ has two structures as a $ U_k $-module, namely $ M^*$ and 
$ M^t $ defined as 
$$ M^*: \,\,\,\,\,\,\,\,\,\,\, (uf)(m) := f(S(u)m) $$
$$ M^t: \,\,\,\,\,\,\,\,\,\,\, (uf)(m) := f(S^{-1}(u)m) $$
where $ S $ is the antipodal map of the Hopf algebra $ U_k $. When $ M $ 
is of finite dimension we have the isomorphisms 
$ M^{*t} \cong M^{t*} \cong M $.

\subsection{}
Using this we can deduce the following Lemma
\begin{lemma}\label{Du} Let $ N $ be a finite 
dimensional $U_k( \fb )$-module. Then 
$$ ( {\cal D } N )^* \cong \mbox{H}_k^0 (N^*) $$
\end{lemma}
\begin{pf*}{Proof} Let $ \Phi \in \mbox{Hom}_{U_k}(\, ( {\cal D } N )^*,
\mbox{H}_k^0 (N^*) \, ) $ be the map corresponding to $ \sigma \in
\mbox{Hom}_{U_k}(\, N, {\cal D } N \, ) $ under the isomorphisms 
$$ \mbox{Hom}_{U_k(\fb)}(\, N, {\cal D } N \, ) \cong
\mbox{Hom}_{U_k(\fb)}(\, ( {\cal D } N)^*, N^* \, )\cong
\mbox{Hom}_{U_k(\fg)}(\, ( {\cal D } N)^*,\mbox{H}_k^0( N^*) \, ) $$
The second isomorphism was Frobenius reciprocity for $ \mbox{H}_k^0 $. 

Let 
$ \Psi \in \mbox{Hom}_{U_k(\fg)}(\, ( \mbox{H}_k^0( N^*),( {\cal D } N)^* \, )$
be the map corresponding to 
$ Ev \in \mbox{Hom}_{U_k(\fb)}(\, \mbox{H}_k^0( N^*), N^* \, )$ under the 
isomorphisms 
$$ \mbox{Hom}_{U_k}(\,\mbox{H}_k^0( N^*),( {\cal D} N)^* \, ) \cong 
\mbox{Hom}_{U_k}(\, {\cal D} N,\mbox{H}_k^0( N^*)^t \, ) \cong $$ 
$$ \mbox{Hom}_{U_k(\fb)}(\, N,\mbox{H}_k^0( N^*)^t \, ) \cong
\mbox{Hom}_{U_k(\fb)}(\, \mbox{H}_k^0( N^*), N^*  \, ) $$
Here the second isomorphism was Frobenius reciprocity for $ \cal D $; it 
can be applied since $ \mbox{H}_k^0( N^*) $ according to [AW] is finite 
dimensional. We can describe $ \Phi $ and $ \Psi $ in the following way
$$ \Phi: \,\,\,\, f \mapsto [ u \mapsto ( n \mapsto f(u(\sigma n )))] $$
$$ \Psi: \,\,\,\, f \mapsto [ u \sigma (n) \mapsto f(u)(n) ] $$
Using these descriptions one checks that $ \Phi \circ \Psi = Id $ and 
$ \Psi \circ \Phi = Id $. 
\end{pf*}

\subsection{} The functor $ \cal D $ is coinvariant and right exact so we 
would like to introduce its derived functors. However, as the category of 
finite dimensional $ U_k $-modules does not have enough projectives one cannot
proceed as normal when defining these. To overcome this obstacle, Joseph 
proposes to use objects of the form $ E \otimes \la $, where $ E $ is 
finite dimensional and $ \la $ is dominant, as substitutes for projectives
[J]. We shall show that this definition also makes sense in our context.
Our methods were here inspired by [P].

\begin{lemma} Let $ M $ be a finite dimensional $ U_k ( \fb ) $-module. Then
there exists a $ \la \in P^+ $ and a finite dimensional $ U_k $-module $ E $ 
such that $ M $ is a quotient of $ E \otimes \la $. 
\end{lemma}

\begin{pf*}{Proof} In [AW] the following ampleness property of $ \mbox{H}^0 $
is shown
\begin{equation} \label{H>>}
\mbox{H}^i ( \la ) = 0 \mbox{ for } i > 0 \mbox{ and } \la << 0
\end{equation}
( We are here inducing from positive Borel groups, hence dominant is 
replaced by antidominant ). As $ M $ is finite dimensional, we can use the 
above to find a $ \la $ such that $ P( M \otimes - \la )  \subset P^- $ and 
such that 
\begin{equation}\label{H1}
\mbox{H}^1 ( \la ) = 0 \, \, \, \, \, \forall \mu \in P ( M  \otimes -\la )
\end{equation}
We now proceed by induction on the cardinality of $  P( M \otimes - \la ) $.
Choose $ \nu $ maximal in $ P( M \otimes - \la ) $; then $ \nu \subset 
M \otimes - \la $ as $ U_k ( \fb ) $-modules. From this we obtain a 
commutative diagram as follows 
$$
\begin{matrix}

0 & \to & \mbox{H}^0_k(\nu) & \to & \mbox{H}^0_k(M \otimes - \la ) & \to 
&  \mbox{H}^0_k((M \otimes - \la)/\nu ) & \to & 0 \\
& & \downarrow \text{\tiny \em Ev} & &  \downarrow \text{\tiny \em Ev} & & 
\downarrow  \text{\tiny \em Ev} & & \\
0 & \to & \nu & \to & M \otimes - \la  & \to 
&  ( M \otimes - \la)/\nu  & \to & 0
\end{matrix}
$$
The rows are both exact, the first one by ~\ref{H1}, the second by 
construction. The first vertical map is surjective by definition and the 
third is surjective by induction hypothesis. But then also the second vertical
map must be surjective and we are done.
\end{pf*}

\subsection{}\label{reso}

We can now in the usual way construct resolutions of all finite dimensional 
$ U_k $-modules $ N $; these will be on the form 
$ (E \otimes \la )^{\bullet} \sur N $ 
and in general infinite. We can furthermore assume that in 
$ (E \otimes \la )^{\bullet} \sur N $ all $ - \la $ satisfy 
$ \mbox{H}^0_k( \la ) 
= 0 $ for $ i > 0 $, this is possible by ~\ref{H>>}. 

\begin{lemma} \label{Dual} Let $ (E \otimes \la )^{\bullet} \sur N $ be a 
resolution
like the above. Then the cohomology $ {\cal D}^{\bullet} $ of the complex
$ {\cal D } ( (E \otimes \la )^{\bullet} \sur 0 $ is independent of the 
choice of resolution. Furthermore there is an isomorphism of $ U_k $-modules 
$$ {\cal D }^i (N)^* \cong  \text{H}^i_k( N^* ) $$
\end{lemma}

\begin{pf*}{Proof} Dualizing the resolution  
$ (E \otimes \la )^{\bullet} \sur N $ we get the resolution 
$ N^* \hookrightarrow (E^* \otimes - \la )^{\bullet}  $. It is acyclic for 
$ \mbox{H}^0_k $ because the tensor identity gives that 
$$ \mbox{H}^i_k ( E^* \otimes  \la ) \cong E^* \otimes 
\mbox{H}^i_k ( \la ) \cong 0 $$
But then $ \mbox{H}^i_k ( N^* )$ is the $i$'th cohomology of $ 0
\hookrightarrow  \mbox{H}^0_k ( E^* \otimes \la )^{\bullet} $, which by 
Lemma~\ref{Du} is the $ i$'th cohomology of 
$ 0 \hookrightarrow ( {\cal D} ( E \otimes \la) )^{* \bullet}$.
The Lemma is proved.
\end{pf*}

\section{The vanishing theorem}

\subsection{}

Assume that we are in the rank one case. i.e. $ \fg = {\fsl}_2 $. We then
have the following well known results

\begin{lemma} \label{D>>} \begin{enumerate} 
\item Let $ \la \geq -1 $. Then $ {\cal D}^j \la = 0 $ for $ j > 0 $, while 
$ {\cal D } \la $ has dimension $ \la + 1 $, a basis being 
$ f^{ (k)} \otimes v_{\la} $. The action of $ U_q ( \fg ) $ is as in 2.2
\item Let $ \la $ be as above and let $ Q $ be the $ U_k ( \fb ) $-module 
$ {\cal D } \la / \la $. Then $ {\cal D } Q = 0 $ and $ {\cal D }^j Q = 0 $
for $ j > 0 $. 
\end{enumerate}
\end{lemma}
\begin{pf*}{Proof} (1) follows from the corresponding Proposition 4.2 in 
[APW] and the preceding Lemma. As for (2) consider the long exact sequence 
of $ U_k ( \fg ) $-modules
$$ \longrightarrow {\cal D }^1 {\cal D } \la 
\longrightarrow {\cal D }^1 Q 
\longrightarrow {\cal D } \la 
\stackrel{g}{\longrightarrow }  {\cal D } \la 
\longrightarrow {\cal D } Q \longrightarrow 0 $$
which arises from the application of $ \cal D $ to the sequence defining 
$ Q $. The definition of derived functors in the last section gives that 
$ {\cal D }^j {\cal D } \la $ for $ j > 0 $; thus the $ j > 1 $ case of 
(2) follows by combining with (1). Now, $ g $ must be a nonzero scalar times 
the identity map on $ {\cal D } \la $ because it is the map corresponding 
to $ \sigma \not= 0 $ under Frobenius reciprocity and 
$$ \Hom_{U_k ( \fg )} ( {\cal D } \la, {\cal D } \la ) \cong  
\Hom_{U_k ( \fb )} ( \la, {\cal D } \la )  \cong k $$
Hence $ { \cal D}^1 Q $ as well as $ { \cal D} Q $ must be zero. The Lemma
is proved.
\end{pf*}

\subsection{} We still consider the rank one case. For a copy of $ {\fsl}_2 $
in $ U_k $ corresponding to the simple reflection $ s $ we denote the 
corresponding Joseph induction $ {\cal D}_s $. 

\begin{theorem} \label{D>>>}\begin{enumerate}
\item Let $ V_w^k ( \la ) $ be as in ~\ref{VW}. Then 
$ {\cal D}_s^j V_w^k ( \la ) = 0 $ for $ j > 0 $.
\item Let $ sw < w $. Then $ {\cal D}_s  V_{sw}^k ( \la ) = 
 V_{w}^k ( \la )  $
\end{enumerate}
\end{theorem}

\begin{pf*} {Proof} ad (1). Let $ i \in I $ be the index corresponding to 
$ s $ and consider the $ W $-filtration of $ V_w^k ( \la ) $ from 
Lemma~\ref{Fi} with respect to this $ i $. The quotients $ Q $ of this all 
satisfy $ {\cal D}_s^j Q = O $ for $ j > 0 $: if $ Q $ is a dominant line 
this is because of Lemma~\ref{D>>} and if $ Q $ is the restriction of an 
$ {\fsl}_2 $-module it follows from the definition of $ {\cal D }^j $ in 
the previous section. Using induction on the filtration length we conclude 
that $ {\cal D}_s^j V_w^k ( \la ) = 0 $ for $ j > 0 $.

\medskip
\noindent
 ad (2). Let the $ U_k( {\fb}_i ) $-module $ R $ be defined as such that 
the following is exact
\begin{equation} \label{Se}
0 \longrightarrow V_{sw}^k ( \la ) 
\longrightarrow V_{w}^k ( \la ) \longrightarrow R 
 \longrightarrow 0
\end{equation}

We know from ~\ref{VW} that $ V_{w}^k ( \la ) $ has a $ k $-basis on the 
form $ \{ G_{\la} ( b ) | b \in B_w ( \la ) \} $; hence $ R $ has a basis 
on the form $ \{ G_{\la} ( b ) | b \in B_w ( \la ) \backslash B_{sw} ( \la )
 \} $. In the $ W $-filtration of $ R $, the quotient between 
the $ l $'th and $ (l+1)$'th term hence has basis 
$ \{ G_{\la} (b) | b \in I^l(B_w ( \la )) \backslash I^l(B_{sw}(\la)) \} $

\medskip
We saw in the proof of Lemma~\ref{Fi} that the $ I^l(B_y ( \la ) ) $ satisfy 
string properties like ~\ref{St}. If $ S$ is an $i$-string and $ S \cap 
I^l ( B_{sw} ( \la ) ) \not= \emptyset  $, then by ~\ref{St} we 
have $ S \subset 
I^l ( B_{sw} ( \la ) ) $. Now $ V_w^k(\la) $ is a $ U_k ( \fg_i ) $-module, 
so $ I^l ( B_w ( \la ) ) $ is the disjoint union of $ i $-strings; and if 
$ S \subset I^l ( B_w ( \la ) ) $ is such a string then from ~\ref{St} 
we get $ S \cap I^l ( B_{sw} ( \la ) ) \not= \emptyset $. Putting these things 
together we see that $ I^l(B_w ( \la )) \backslash I^l(B_{sw}(\la)) $ is the 
union of $ i $-strings with the highest weight vector omitted; this must 
furthermore be a of $ \fsl_2 $-weight $ l $. From the formula ~\ref{V2} we
then see that the span of the global basis elements of $ I^l ( B_w ( \la ))
\setminus B_{sw}  ( \la ) $. form a $ U_k ( \fb_i ) $-module of type 
$ { \cal D \la } / \la $. And then an induction on the filtration length, 
the induction start being provided by (2) of Lemma~\ref{D>>} proves that 
$ {\cal D }^1_s R = { \cal D }_s R = 0 $. 

\medskip

Now, $ V_w^{k}( \la ) $ is a $ U_k ( \fg_i ) $-module so we have 
$ { \cal D }^1 R = { \cal D }_s R = 0 $ we may finish the proof of 
(2) by applying $ \cal D $ to ~\ref{D>>>} (1).
\end{pf*}

In Theorem~\ref{D>>>} the induction $ {\cal D }_s $ was induction from a
rank one Borel subgroup. However, $ V_w^k(\la)$ has a module structure for 
the full Borel group $ U_k(\fb) $. Denote by $ U_k ( \fg_i ) $ the minimal 
parabolic quantum group generated by this and the $ f^{(n)}_i $'s and by 
$ {\cal D }^{\prime}_s $ the induction from $ U_k( \fb ) $ to $ U^{\prime}_k
( {\fg}_i ) $. As $k$-vector spaces and $ U_k ( \fg_i ) $-modules 
$  {\cal D }^{\prime}_s V_w^k ( \la ) $ and $ {\cal D }_s V_w^k ( \la ) $
are isomorphic, namely both equal to the largest $ f_i $-finite quotient
of $$ U_k ( \fg_i ) \otimes_{ U_k( \fb_i ) } V_{sw}^k ( \la ) $$
We have the following Lemma

\begin{lemma}\label{U'} There is an isomorphism of 
$ U^{\prime}_k( \fg_i ) $-modules
$$   {\cal D }^{\prime}_s V_{sw}^k ( \la ) \cong V_w^k ( \la ) $$
\end{lemma} 

\begin{pf*}{Proof} There is a commutative diagram 
$$ \begin{matrix}
V_{sw}^k ( \la) & \stackrel{id}{\longrightarrow} & V_{sw}^k ( \la) \\
\text{ \tiny{$ \sigma $ }} \downarrow & &  \downarrow \text{ \em \tiny i } \\
{\cal D }^{\prime}_s V_{sw}^k ( \la) & \stackrel{\varphi}{\longrightarrow} 
& V_{w}^k ( \la)
\end{matrix} $$
where $ i $ is the inclusion map, $ \sigma $ is the canonical map -- these 
are $ U^{\prime}_k ( \fg_i ) $-linear -- and $ \varphi $ is the 
$ U_k ( \fg_i ) $-linear map obtained from Theorem~\ref{D>>>} (2) together 
with the remarks in the beginning of this section.

\medskip

Any element of $ { \cal D }^{\prime}_s V_{sw}^k ( \la )  $ can be written 
as a linear combination of elements on the form $ u \sigma ( v ) $ where 
$ u \in U_k^- ( \fg_i ), v \in V^k_{sw} ( \la ) $. We must show that 
$ \varphi $ commutes with $ e_j $ for $ j \not= i $; it suffices to do so 
for $ u \sigma( v ) $. On the one hand we have 
$$ \varphi( e_j ( u \sigma ( v))) = \varphi ( u e_j \sigma ( v)) = 
\varphi ( u \sigma ( e_j v )) = u \varphi ( \sigma ( e_j v )) = 
u i ( e_j v ) = e_j ui ( v) $$
where we used the commutativity of the diagram and of $ e_j $ and $ f_i $.
On the other hand 
$$ e_j \varphi ( u \sigma (v)) = e_j u \varphi ( \sigma ( v)) = 
e_j u i ( v) $$
We see that the two sides are equal.
\end{pf*}

\subsection{} We omit from now on the primes on the $ \cal D $, inductions 
will be from the full Borel subalgebra. For $ w_0 = s_{1_1} \, s_{1_2}\,  
\ldots \, s_{1_n} $ a reduced expression of the longest element of $ W $ 
we define $ {\cal D }_{w_0} $ as $ {\cal D }_{i_1} \, 
{\cal D }_{i_2} \, \ldots {\cal D }_{i_n} $, ( apriori this may depend on 
the chosen expression ). The Weyl module $ V_k ( \la ) $ is by definition 
$ V_{ \DZ } ( \la) \otimes k $. We now obtain the following Theorem
 
\begin{theorem}\label{Iso} For $ \la \in P^+ $ there are isomorphisms of 
$ U_k( \fb )$-modules 
$$ { \cal D}_{w_0} \la \cong  { \cal D} \la \cong V_k ( \la ) $$ 
\end{theorem}

\begin{pf*} {Proof} As $ V_k (\la ) $ has finite dimension there is by 
definition of $ \cal D $ a surjection 
$ \varphi: {\cal D } \la \sur V_k( \la ) $. Now, succesive applications of 
Theorem~\ref{D>>>} and Lemma~\ref{U'} give the isomorphism 
$ { \cal D}_{w_0} \la \cong V_k (\la )$. And applying
Frobenius 
reciprocity succesively, the canonical map $ \sigma: \la \rightarrow 
{\cal D } \la $ induces a $ U_k( \fb ) $-linear map 
$ \psi: { \cal D}_{w_0} \la \cong V_k (\la ) \rightarrow {\cal D } \la $.
But any $  U_k( \fb ) $-linear map must also be a $  U_k( \fg ) $-linear map
as one sees from Frobenius reciprocity: $ V_k ( \la ) $ and $ {\cal D } \la $ 
are both $  U_k( \fg ) $-modules. Thus the composition $ \varphi \circ 
\psi $ is a $ U_k ( \fg ) $-linear endomorphism of $ V_k ( \la ) $ and one 
checks that it is nonzero on the $ \la $'th weight space. But $ V_k(\la ) $
is a highest weight module, hence $  \varphi \circ \psi $ must be a 
nonzero scalar times the identity. This shows that $ {\cal D } \la 
\cong V_k( \la ) \oplus M $ for $ M $ some $ U_k ( \fg ) $-module. Now, 
$ {\cal D }\la $ is indecomposable being a highest weight module as well, and 
we get a contradiction unless $ M = 0 $. The Theorem is proved. 
\end{pf*}

\subsection{}

We can now prove our main Theorem.
\begin{theorem} ( Kempf vanishing ). Let $ \la \in P^- $. Then $ \text{H}_k^i
( \la ) = 0 $ for $ i > 0 $.
\end{theorem}

\begin{pf*} {Proof} By Lemma~\ref{Dual} the Theorem is equivalent to 
$ {\cal D }^i ( -\la ) = 0 $. Let $ ( E \otimes \nu )^{\bullet}  \rightarrow
- \la \rightarrow 0 $ be a resolution of $ - \la $ as in Section~\ref{reso}.
Then the complex $ ( {\cal D}_{s_n} ( E \otimes \nu ))^{\bullet} 
\rightarrow {\cal D}_{s_n} ( -\la ) \rightarrow 0 $ is exact by 
Lemma~\ref{D>>} (1). It is also acyclic for $ {\cal D }_{s_{n-1} } $ because
by the tensor identity and Theorem~\ref{D>>>} we have 
$$ {\cal D}^i_{s_{n-1}} ( {\cal D}_{s_{n}} ( E \otimes \nu )) \cong 
E \otimes {\cal D}^i_{s_{n-1}} ( {\cal D}_{s_{n}} ( \nu )) \cong 
E \otimes {\cal D}^i_{s_{n-1}} ( {\cal D}_{s_{n}} ( V^k_{s_n}(\nu ))
\cong 0 $$
Thus the application of $ {\cal D }_{s_{n-1}} $ to 
$ ( {\cal D}_{s_{n}} ( E \otimes \nu ))^{\bullet} \rightarrow 
  {\cal D}_{s_{n}} ( -\la ) \rightarrow 0 $ gives a complex that evaluates 
$ {\cal D}^i_{s_{n-1}} ( {\cal D}_{s_{n}} ( -\la )) $. But we know from
Theorem~\ref{D>>>} that these cohomology groups are zero so the complex 
$ {\cal D}_{s_{n-1}} ( {\cal D}_{s_{n}} ( E \otimes \nu ))^{\bullet}
\rightarrow 
{\cal D}_{s_{n-1}} ( {\cal D}^i_{s_{n}} ( - \la  )) \rightarrow 0 $ is exact.
Furthermore, the argument from before shows that it is acyclic for 
$ {\cal D}_{s_{n-2}} $. Continuing we eventually reach the sequence 
$ ({ \cal D }_{w_0} ( E \otimes \nu ))^{\bullet} \rightarrow 
{\cal D}_{w_0} ( - \la ) \rightarrow 0 $ which thus is exact. By 
Theorem ~\ref{Iso} it is isomorphic to 
$ ({ \cal D } ( E \otimes \nu ))^{\bullet} \rightarrow 
{\cal D} ( - \la ) \rightarrow 0 $ and we are done.
\end{pf*}

\subsection{}

We have a couple of corollaries to this Theorem.

\begin{corollary} ( Demazure vanishing ). For $ \la \in P^- $ we have 
$$ { \cal D }^i_{s_n} ( { \cal D }_{s_{n-1}} { \cal D }_{s_{n-2}} \ldots
{ \cal D }_{s_{1}} (\la)) \cong 0 $$
\end{corollary}
\begin{pf*} {Proof} This is contained in the last proof.
\end{pf*}

\begin{corollary} The modular Kempf and Demazure vanishing theorems.
\end{corollary}
\begin{pf*} {Proof} The result follows from the Theorem and the Corollary 
by specializing $ q $ to $ 1 $; there are base change theorems controlling 
this. 
\end{pf*}

\begin{corollary} Demazure's character formula in terms of the $ \mbox{H}^i $.
\end{corollary}
\begin{pf*} {Proof} The classical proof carries over, see [A1] 
\end{pf*}

\end{document}